\documentclass[11pt]{article}
\usepackage{tabularx}
\usepackage{graphicx}
\usepackage{fullpage}
\usepackage{hyperref}
\usepackage{natbib}
\usepackage{delarray}
\usepackage{times}
\usepackage{subfigure}
\usepackage{color}
\usepackage{amssymb}
\usepackage{amsmath}
\usepackage{tikz}
\usepackage{xspace}
\usepackage{enumerate}
\usepackage{algorithm}
\usepackage{algorithmic}
\usepackage{multirow}
\usepackage{nameref}
\usepackage{booktabs}
\usepackage{mdframed}
\usepackage{fmtcount}
\usepackage{footnote}
\usepackage{bm}

\newcommand{\eat}[1]{}
\newcommand{\topic}[1]{\vspace{5pt}\noindent{{\bf #1:}}}

\newcommand{\R}{{\mathbb R}}

\newcommand{\red}[1]{\textcolor{red}{#1}}
\newcommand{\blue}[1]{\textcolor{blue}{#1}}
\newcommand{\inner}[2]{\langle #1, #2 \rangle}

\newcommand{\n}[1]{\| #1 \|}

\newcommand{\tx}{\widetilde{x}}

\newcommand{\bx}{\bar{x}}

\newcommand{\DVZ}{Dadush-V{\'e}gh-Zambelli }
\newcommand{\tbp}{\mathsf{T_{BP}}}
\newcommand{\tma}{\mathsf{T_{MA}}}
\newcommand{\ppdp}{\mathsf{PPDP}}
\newcommand{\1}{\mathbf{1}}
\newcommand{\ths}{\mathrm{threshold}}
\newcommand{\case}{\mathrm{case}}
\newcommand{\diag}{\mathrm{diag}}
\newcommand{\rank}{\mathrm{rank}}
\newcommand{\bound}{\mathrm{bound}}
\newcommand{\rt}{\textbf{return}}

                {\hspace*{\fill}$\Box$\par}
        {\hspace*{\fill}$\Box$\par\vspace{4mm}}
        {\hspace*{\fill}$\Box$\par}

\newcommand*{\QED}{\hfill\ensuremath{\square}}

\newtheorem{theorem}{Theorem}
\newtheorem{lemma}{Lemma}

\hypersetup{
    colorlinks=true,       
    linkcolor=black,    
    citecolor=black,        
    filecolor=magenta,     
    urlcolor=blue
}
\begin{document}

\begin{titlepage}

\title{A Fast Polynomial-time Primal-Dual Projection Algorithm for Linear Programming}
\author{Zhize Li \\ IIIS, Tsinghua University \\ zz-li14@mails.tsinghua.edu.cn
        \and Wei Zhang \\ TLI, National University of Singapore \\ lindelfeel@gmail.com
        \and Kees Roos \\ EEMCS, Delft University of Technology \\ c.roos@tudelft.nl}

\date{}
\clearpage\maketitle
\thispagestyle{empty}

\begin{abstract}
Traditionally, there are several polynomial algorithms for linear programming including the ellipsoid method, the interior point method and other variants.
Recently, Chubanov \citep{chubanov2015polynomial} proposed a projection and rescaling algorithm, which has become a potentially \emph{practical} class of polynomial algorithms for linear feasibility problems and also for the general linear programming.
However, the Chubanov-type algorithms usually perform much better on the infeasible instances than on the feasible instances in practice.
To explain this phenomenon, we derive a new theoretical complexity bound for the infeasible instances based on the condition number,
which shows that algorithms can indeed run much faster on infeasible instances in certain situations.
In order to speed up the feasible instances, we propose a \emph{Polynomial-time Primal-Dual Projection} algorithm (called $\ppdp$) by explicitly developing the dual algorithm.
The numerical results validate that our $\ppdp$ algorithm achieves a quite balanced performance between feasible and infeasible instances, and its performance is remarkably better than previous algorithms.
\end{abstract}

\end{titlepage}

\section{Introduction}
\label{sec:intro}

Linear programming is a fundamental problem in many areas, such as operations research, network, machine learning, business analysis and finance \citep{von1947theory,dantzig1963linear,luenberger1984linear,boyd2004convex}.
In this paper, we consider the \emph{maximum support} of the linear feasibility problem
\begin{equation}
\label{feasibility:primal}
\begin{aligned}
\mathrm{find} ~~&x \in \R^n \\
\mathrm{subject~to}~~ &Ax=0, ~x\geq 0,~ x\neq 0,
\end{aligned}
\end{equation}
with its dual problem
\begin{equation}
\label{feasibility:dual}
\begin{aligned}
\mathrm{find} ~~&u \in \R^m \\
\mathrm{subject~to}~~ &A^T u> 0,
\end{aligned}
\end{equation}
where $A\in \R^{m\times n}$ is an integer (or rational) matrix and $\rank(A) = m$.
The \emph{maximum support} means that the set of positive coordinates of the returned solution of (\ref{feasibility:primal}) should be inclusion-wise maximum.
Actually, for the solution $\hat{x}$ returned by our algorithm, any coordinate $\hat{x}_i=0$ if and only if this coordinate equals to 0 for all feasible solutions of (\ref{feasibility:primal}).
Thus, our algorithm can be directly used to test the feasibility of the general linear system $Ax=b,x\geq 0$ with the same time complexity, i.e., given the maximum support solution $(\bx,\bx')$ to the system $Ax-bx'=0, (x,x')\geq 0$, if $\bx'>0$ then the original problem  $Ax=b,x\geq0$ has a solution $\tx = \bx/\bx'$, otherwise it is infeasible.

There are many (polynomial-time) algorithms for solving linear programming problems, e.g., \citep{karmarkar1984new}, \citep{wright1997primal} and \citep{renegar1988polynomial}.
Recently, \cite{chubanov2015polynomial} proposed a polynomial-time projection and rescaling algorithm for solving problem \eqref{feasibility:primal}.
Due to its simplicity and efficiency, this kind of algorithms has become a potentially \emph{practical} class of polynomial algorithms. See e.g., \citep{dadush2016rescaled}, \citep{roos2018improved} and \citep{pena2018computational}.

Chubanov's algorithm \citep{chubanov2015polynomial} and its variants typically consist of two procedures. The key part is \emph{basic procedure} (BP) and the other part is \emph{main algorithm} (MA).
The BP returns one of the following three results:
\begin{enumerate}[(i)]
  \item a feasible solution of \eqref{feasibility:primal};
  \item a feasible solution of the dual problem \eqref{feasibility:dual};
  \item a cut for the feasible region of \eqref{feasibility:primal}.
\end{enumerate}
Note that exactly one of \eqref{feasibility:primal} and \eqref{feasibility:dual} is feasible according to Farkas' lemma. Thus \eqref{feasibility:primal} is infeasible if BP returns (ii).
If BP returns (iii), the other procedure MA rescales the matrix $A$ by using this cut and call BP again on the rescaled matrix $A$.
According to \citep{khachian1979polynomial} which gives a positive lower bound on the entries of a solution of a linear system, after a certain number of rescalings, one can conclude that there is no feasible solution for \eqref{feasibility:primal}. So the number of rescaling operations can be bounded, i.e., the number of MA calls can be bounded.
Consequently, the algorithm can terminate in finite time no matter whether problem \eqref{feasibility:primal} is feasible or infeasible.

To be more precise, we quantify the time complexity.
The total time complexity of these Chubanov-type algorithms are typically $O(\tma*\tbp)$,
where $\tma$ denotes the number of MA calls (rescaling operations), and $\tbp$ denotes the time required by the basic procedure BP.
According to the classic lower bound \citep{khachian1979polynomial}, $\tma$ can simply be bounded by $O(nL)$ for these Chubanov-type algorithms, where $L$ denotes the bit size of $A$.
However, $\tbp$ is the most important and tricky part.
Theoretical results and practical performances vary for different BP procedures.
The typical BP procedures include the perceptron method, von Neumann's method, and their variants (see e.g., \citep{dantzig1992varepsilon,dunagan2008simple,dadush2016rescaled,pena2017solving}).
We review more details of these BP in Section \ref{sec:bp}.
Usually, $\tbp$ equals to $O(n^4)$ or $O(n^3m)$ in these BP procedures.
In this work, we improve $\tbp$ by a factor of $\sqrt{n}$ if (\ref{feasibility:primal}) or (\ref{feasibility:dual}) is well-conditioned (measured by \eqref{eq:rho}), but in the worse case, $\tbp$ still equals to $O(n^3m)$ in our algorithm.

\topic{Our Motivation}
In practice, these Chubanov-type projection and rescaling algorithms usually run much faster on the primal infeasible instances (i.e., \eqref{feasibility:primal} is infeasible) than on the primal feasible instances (i.e., dual infeasible) no matter what basic procedure (von Neumann, perceptron or their variants) we use (also see Table \ref{table:comparison} in Section \ref{sec:exp}).
In this paper, we try to explain this phenomenon theoretically. Moreover, we try to provide a new algorithm to address this issue.

\topic{Our Contribution} Concretely, we make the following technical contributions:
\begin{enumerate}
  \item First, for the theoretical explanation, we provide Lemma \ref{lm:infeasible} which shows that
  the time complexity $\tbp$ can be $O(n^{2.5}m)$ rather than $O(n^3m)$ (see Lemma \ref{lm:tbp}) in certain situations if \eqref{feasibility:primal} is infeasible.
  This gives an explanation of why these Chubanov-type algorithms usually run much faster if \eqref{feasibility:primal} is infeasible.
  \item Then,  we explicitly develop the \emph{dual algorithm} (see Section \ref{sec:dual}) to improve the performance when \eqref{feasibility:primal} is feasible.
      Our dual algorithm is the first algorithm which rescales the row space of $A$ in MA (see Table \ref{tab:comp}).
      As a result, we provide a similar Lemma \ref{lm:infeasibled} which shows that the time complexity $\tbp$ of our dual algorithm can be $O(n^{2.5}m)$ rather than $O(n^3m)$ in certain situations if \eqref{feasibility:dual} is infeasible (i.e. \eqref{feasibility:primal} is feasible).

        Naturally, we obtain a new fast polynomial primal-dual projection algorithm (called $\ppdp$) by integrating our primal algorithm (which runs faster on the primal infeasible instances) and our dual algorithm (which runs faster on the primal feasible instances). See Section \ref{sec:ppdp}.
  \item Finally, the numerical results validate that our primal-dual $\ppdp$ algorithm is quite balanced between feasible and infeasible instances, and it runs significantly faster than other algorithms (see Table \ref{table:comparison} in Section \ref{sec:exp}).
\end{enumerate}

\topic{Remark}
Our algorithms are based on \DVZ algorithm \citep{dadush2016rescaled} and the improvements of Roos's algorithm \citep{roos2018improved} (see Section \ref{sec:relatealg} and Table \ref{tab:comp}).
Besides, we introduce a new step-size term $c$ for practical consideration (see Line 13 and 14 of Algorithm \ref{alg:bp} and \ref{code:dual Dadush}).
For the maximum support problem \eqref{feasibility:primal}, $\tbp=O(n^4)$ for Chubanov's algorithm and Roos's algorithm, and  $\tbp=O(n^3m)$ for \DVZ algorithm.
Note that in the worst case $\tbp=O(n^3m)$ for our algorithms, but it can be improved by a factor of $\sqrt{n}$ in certain situations.
Recall that $\tma=O(nL)$ for these Chubanov-type algorithms as we discussed before.
Thus the time complexity of our algorithms (in the worst case) match the result of \DVZ algorithm \citep{dadush2016rescaled}, i.e., $O(\tma*\tbp)=O(n^4mL)$ (see our Theorems \ref{thm:primal}--\ref{thm:primal_dual}).
However, we point out that the total time complexity of Chubanov’s algorithm and Roos's algorithm are $O(n^4L)$, and hence is faster than ours. They speed up it from $O(n^5L)$ to $O(n^4L)$ by using an amortized analysis while we currently do not use. We leave this speedup as a future work.

\topic{Organization}
In Section \ref{sec:pre}, we introduce some useful notations and review some related algorithms.
The details and results for our primal algorithm and dual algorithm are provided in Section \ref{sec:primal} and Section \ref{sec:dual}, respectively.
Then, in Section \ref{sec:ppdp}, we propose the efficient primal-dual $\ppdp$ algorithm.
Finally, we conduct the numerical experiments in Section \ref{sec:exp} and include a brief conclusion in Section \ref{sec:con}.

\section{Preliminaries}
\label{sec:pre}

In this section, we first review some classic basic procedures and then introduce some notations to review some related algorithms at the end of this section.
\subsection{Classic Basic Procedures} \label{sec:bp}
Recall that BP returns one of the following three results:
(i) a feasible solution of \eqref{feasibility:primal};
(ii) a feasible solution of the dual problem \eqref{feasibility:dual};
(iii) a cut for the feasible region of \eqref{feasibility:primal}.
Here we focus on the first two outputs, the last one is controlled by an upper bound lemma (similar to Lemma \ref{lm:roosbound}).

Letting $y=Ax$, when solving \eqref{feasibility:dual}, we know that there is at least an index $k$
such that $a_k^Ty\leq 0$, where $a_k$ is the $k$th-column of $A$ (otherwise $y$ is already a feasible solution for \eqref{feasibility:dual}).
On the other hand, to solve \eqref{feasibility:primal}, we want to minimize $\n{y}$. The goal is to let $y$ go to 0 (in which case $x$ is a feasible solution for \eqref{feasibility:primal}).
We review some classic update methods as follows:

\topic{von Neumann's algorithm}
In each iteration, find an index $k$ such that $a_k^Ty\leq 0$, and then update $x$ and $y$ as
\begin{equation}\label{eq:update}
\red{y'}=\alpha y+\beta a_k, \quad x'=\alpha x+\beta e_k \quad (\mathrm{note~ that~}y=Ax \mathrm{~and~} \red{y'}=Ax'),
\end{equation}
where $\alpha,\beta>0$ are chosen such that $\|\red{y'}\|$ is smallest and $\alpha+\beta=1$ \citep{dantzig1992varepsilon}.
\begin{figure}[!h]
\centering\includegraphics[scale=0.8]{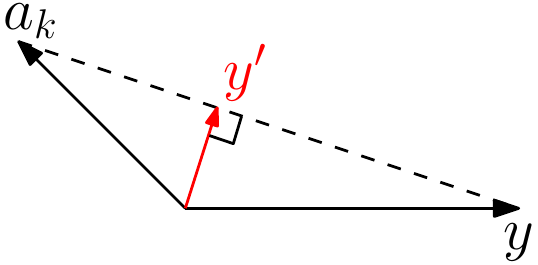}
\end{figure}

\topic{Perceptron}
Choose $\alpha=\beta=1$ in \eqref{eq:update} at every iteration. See e.g. \citep{rosenblatt1957perceptron,novikoff1963convergence}.

\topic{Dunagan-Vempala}
Fix $\alpha=1$ and choose $\beta$ to minimize $\n{\red{y'}}$ \citep{dunagan2008simple}.

\subsection{Notations}
Before reviewing the related algorithms (in the following Section \ref{sec:relatealg}), we need to define/recall some useful notations.
We use $P_A$ and $Q_A$ to denote the projections of $\R^n$ onto the null space ($\mathcal{N}_A$) and row space ($\mathcal{R}_A$) of the $m\times n$ matrix $A$, respectively:
\begin{equation*}
P_A \triangleq I-A^T (A A^T)^{\dag}A,\qquad Q_A \triangleq A^T (A A^T)^{\dag}A.
\end{equation*}
where $(\cdot)^{\dag}$ denotes the Moore-Penrose pseudoinverse. Particularly, $(A A^T)^{\dag}=(A A^T)^{-1}$ if $\rank(A)=m$.

We further define the following notations:
\begin{equation}
\label{eq:vec}
v=Q_Ay\in \mathcal{R}_A,\quad z=P_Ay\in \mathcal{N}_A,\quad y=v+z\in \R^n.
\end{equation}
Usually, $z$ is used to denote the feasible solution of \eqref{feasibility:primal} and $v$ indicates the feasibility of \eqref{feasibility:dual}.

To analyze case (iii) of BP, we note that \eqref{feasibility:primal} is feasible if and only if the system
\begin{equation}
\label{feasibility:normalized primal}
\begin{aligned}
Ax=0, ~x\in [0,1]^n,~ x\neq 0
\end{aligned}
\end{equation}
is feasible since \eqref{feasibility:primal} is a homogeneous system.
From now on, we will consider problem \eqref{feasibility:normalized primal} instead of \eqref{feasibility:primal}. Similarly, we use a normalized version \eqref{eq:bounded dual} to replace \eqref{feasibility:dual}.
Now, we recall a useful lemma which gives an upper bound for the coordinates of any feasible solution. This upper bound will indicate a cut for case (iii).
\begin{lemma}[\citep{roos2018improved}]
\label{lm:roosbound}
Let $x$ be any feasible solution of \eqref{feasibility:normalized primal}, $y$ and $v$ are defined as in (\ref{eq:vec}), then every non-zero coordinate $v_j$ of $v$ gives rise to an upper bound for $x_j$, according to
\begin{equation*}
x_j\leq \bound_j(y)\triangleq \1^T\Big[\frac{v}{-v_j}\Big]^+,
\end{equation*}
where $x^+\triangleq \max\{0,x\}$ and $\1$ denotes the all-ones vector.
\end{lemma}
This means that we can scale the column $j$ of $A$ by a factor $\bound_j(y)$ to make the feasible solutions of (\ref{feasibility:normalized primal}) closer to the all-ones vector $\1$.
Similarly to $x^+$, we denote $x^-\triangleq -(-x)^+$.
Furthermore, we need the definition of condition number $\rho(Q)$ for a matrix $Q$ \citep{goffin1980relaxation}:
\begin{equation}\label{eq:rho}
\rho(Q)\triangleq \max_{x:\|x\|_2=1}\min_{i} \langle x, \frac{q_i}{\|q_i\|_2} \rangle,
\end{equation}
where $q_i$ is the $i$th-column of $Q$.

\subsection{Related Algorithms}\label{sec:relatealg}
Now, we are able to review some related algorithms for solving \eqref{feasibility:normalized primal} based on the BP procedures introduced in Section \ref{sec:bp}.

\topic{Chubanov's algorithm}
Instead of updating in the original space $y=Ax$, \cite{chubanov2015polynomial} updates in the projection space $z=\blue{P_A}y$, where $P_A = I-A^T (A A^T)^{-1}A$ is a null space projection of $A$.
In each BP iteration, it updates $y$ and $z$ in the same way as von Neumann's update (just replacing $A$ by $P_A$).
Intuitively, BP either finds a feasible solution $x^*$ of \eqref{feasibility:normalized primal} or finds a cut (i.e., an index $j$ such that $x_j^*\leq 1/2$ for any feasible solution $x^*$ of \eqref{feasibility:normalized primal} in $[0,1]^n$).
Then the main algorithm MA rescales the null space of $A$ by dividing the $j$th-column of $A$ by 2.
According to [Khachian, 1979], there is a lower bound for the feasible solutions of \eqref{feasibility:normalized primal}.
Thus the number of rescaling operations can be bounded.
Finally, the algorithm terminates in polynomial-time, where either BP returns a feasible solution or MA claims the infeasibility according to the lower bound.

\topic{Roos's algorithm} \cite{roos2015chubanov,roos2018improved} provided two improvements of Chubanov's algorithm:
\begin{enumerate}
  \item A new cut condition was proposed, which is proved better than the one used by Chubanov.
  \item The BP can use \emph{multiple indices} to update $z$ and $y$ ($z=P_Ay$) in each iteration, e.g., a set of indices satisfying $(P_A)_i^Tz\leq 0$. Recall that von Neumann's update only uses one index $k$ satisfying $a_k^Ty\leq 0$.
\end{enumerate}

\topic{Dadush-V{\'e}gh-Zambelli}
Compared with Chubanov's algorithm, \cite{dadush2016rescaled} used the Dunagan-Vempala update instead of von Neumann's update as its BP, along with Roos' new cut condition. Besides, the updates are performed in the orthogonal space $v=\blue{Q_A}y$, where $Q_A = A^T (A A^T)^{-1}A$ is a row space projection matrix of $A$.
But the rescaling space in MA is the same, i.e., the null space of $A$.

\topic{Comparison}
To demonstrate it clearly, we provide a comparison of our algorithms with other algorithms in Table \ref{tab:comp}. Note that our primal-dual $\ppdp$ algorithm is the integration of our primal algorithm and dual algorithm.

\begin{table}[!ht]
  \caption{Comparison of our algorithms with other algorithms}
  \vspace{1mm}
  \label{tab:comp}
 \centering
 \begin{tabular}{ccccc}
  \toprule
    Algorithms  & Update method & Update space & Rescaling space&
    \#indices\\
  \midrule
  Chubanov's algorithm &  von Neumann & Null space & Null space &One\\
  Roos' algorithm     & von Neumann & Null space & Null space&Multiple\\
  Dadush-V{\'e}gh-Zambelli  & Dunagan-Vempala &Row space & Null space&One\\
  Our primal algorithm& Dunagan-Vempala  & Row space& Null space&Multiple\\
  Our dual  algorithm & Dunagan-Vempala  & Null space& Row space&Multiple\\
  \bottomrule
 \end{tabular}
\end{table}

\newpage
\section{Our Primal Algorithm}
\label{sec:primal}

In this section, we introduce our primal algorithm which consists of the basic procedure BP and the main algorithm MA.
The details of BP and MA are provided in Section \ref{sec:bpprimal} and Section \ref{sec:maprimal} respectively.

\subsection{Basic Procedure (BP)}\label{sec:bpprimal}
Our BP is similar to \citep{dadush2016rescaled} (or \citep{dunagan2008simple}) (see Table \ref{tab:comp}). The details are described in Algorithm \ref{alg:bp}.
The main difference is that we use multiple indices $K$ to update (see Line 9 of Algorithm \ref{alg:bp}) and introduce the step-size $c$ for practical consideration (see Line 13 and 14 of Algorithm \ref{alg:bp}).

\begin{algorithm}[!htb]
	\caption{Basic Procedure for the Primal Problem}
	\label{alg:bp}
	\begin{algorithmic}[1]
		\REQUIRE $Q_A$
        \ENSURE
            $y, z, J, \case.$
        \STATE $r=size(Q_A), \ths=1/2r^{3/2}, c\in (0,2), \case=0$
        \STATE $y=\1/r, v=Q_Ay, z=y-v=P_Ay$
		\WHILE{$\case=0$}
		\IF {$z>0$}
		    \STATE $\case = 1$ ($z$ is primal feasible); \rt
		\ELSIF {$v> 0\; and\; r==n$}
            \STATE $\case = 2$ ($v$ is dual feasible); \rt
		\ELSE
            \STATE find $K=\{k: v_k \leq 0\}$
		    \STATE $q_K=Q_A \sum_{k\in K}e_k$
		    \STATE $\alpha=\inner{\frac{q_K}{\|q_K\|_2}}{v}$
		    \IF {$\alpha \leq -\ths$}
                \STATE $y=y-c(\frac{\alpha}{\|q_K\|_2}\sum_{k\in K} e_k)$
		        \STATE $v=v-c(\frac{\alpha}{\|q_K\|_2} \sum_{k\in K} q_k)$
		    \ELSE
		        \STATE find a nonempty set $J$ such that $J\subseteq \{j: bound_j(y)\leq \frac{1}{2} \}$ (a cut); \rt
		    \ENDIF
		\ENDIF
		\ENDWHILE
	\end{algorithmic}
\end{algorithm}

In the BP (Algorithm \ref{alg:bp}), the norm of the iterated vector $v=Q_A y$ is decreasing, while each coordinate of $y$ is increasing.
Thus, after a certain number of iterations, we will obtain a feasible solution $z=y-v=P_Ay>0$.
Otherwise, it is always possible to find a cut $J$ (Line 16), along with some rescaling operations for the matrix $A$, to make the feasible solutions of (\ref{feasibility:normalized primal}) closer to the all-ones vector.
The cut is guaranteed by the following lemma.
\begin{lemma}\label{lm:cut}
Let $Q_A$ be the projection matrix at a given iteration of BP (Algorithm \ref{alg:bp}).
Suppose that $\alpha=\inner{\frac{q_K}{\|q_K\|_2}}{v}> -\ths$, then the set $J=\{j: \bound_j(y)\leq \frac{1}{2} \}$ is nonempty and every solution $x$ of problem (\ref{feasibility:normalized primal}) satisfies $x_j\leq \frac12$ for all $j\in J$.
\end{lemma}

This lemma is proved with Lemma \ref{lm:roosbound} and we defer the proof to Appendix \ref{app:lmcut}.

For the time complexity of Algorithm \ref{alg:bp}, i.e. $\tbp$, we give the following lemma (the proof is in Appendix \ref{app:tbp}).
\begin{lemma}\label{lm:tbp}
The time complexity of Algorithm \ref{alg:bp} $\tbp=O(n^3m)$. Concretely, it uses at most $O(n^2)$ iterations and each iteration costs at most $O(mn)$ time.
\end{lemma}

Note that Lemma \ref{lm:tbp} holds regardless \eqref{feasibility:normalized primal} is feasible or infeasible. However, as we discussed before, the algorithm usually performs much better on the infeasible instances than on the feasible instances. To explain this phenomenon, we give the following lemma. The proof is deferred to Appendix \ref{app:infeasible}.
\begin{lemma}\label{lm:infeasible}
If \eqref{feasibility:normalized primal} is infeasible,
the time complexity of Algorithm \ref{alg:bp} $\tbp=O(n^2m/\rho(Q_A))$, where $\rho(Q_A)$ is the condition number defined in (\ref{eq:rho}). 
In particular, $\rho(Q_A)$ equals to $1/\sqrt{n}$ under well-condition (e.g., $A$ is an identity matrix), then $\tbp=O(n^{2.5}m)$ if problem \eqref{feasibility:normalized primal} is infeasible.
\end{lemma}

\subsection{Main Algorithm (MA)}\label{sec:maprimal}
The details of our MA are described in Algorithm \ref{alg:ma}.
Particularly, we rescale the null space of $A$ in Line 8.
\begin{algorithm}[!htb]
	\caption{Main Algorithm for the Primal Problem}
	\label{alg:ma}
	\begin{algorithmic}[1]
		\REQUIRE
            $A\in\R^{m\times n}, d=\1, \tau=2^{-L}, \case=0, H=\varnothing$.
		\WHILE{$\case=0$}
            \STATE $Q_A=A^T (A A^T)^{\dag}A$
            \STATE $(y,z,J,\case)\leftarrow$ Basic Procedure for Primal Problem$(Q_A)$
		\IF {$\case=0$}
		    \STATE $d_J=d_J /2$
		
		   \STATE $H=\{i:d_i\leq \tau\}$
		   \STATE $d_H=0$
		   \STATE $A_{J}=A_{J}/2$
		   \STATE $A=A_{\overline{H}}$
		
		\ENDIF
		\ENDWHILE
        \IF{$\case=1$}
            \STATE $d=d_{\overline{H}}$
            \STATE $D=\diag(d)$
            \STATE Define $x$ as $x_{\overline{H}}=Dz, x_H=0$
        \ENDIF
	\end{algorithmic}
\end{algorithm}

Now, we state the complexity of our primal algorithm in the following theorem.
The proof is deferred to Appendix \ref{app:thmprimal}
\begin{theorem}\label{thm:primal}
The time complexity of the primal algorithm is $O(n^4mL)$.
\end{theorem}

\section{Our Dual Algorithm}
\label{sec:dual}
The Chubanov-type algorithms all focus on the primal problem (\ref{feasibility:primal}) (or the normalized version \eqref{feasibility:normalized primal}), i.e.,
their MA always rescale the null space of $A$ (see Table \ref{tab:comp}).
We emphasize that these algorithms usually perform much better on the infeasible instances than on the feasible ones (see our Lemma \ref{lm:infeasible} which gives an explanation).
Now, we want to address this unbalanced issue by providing a dual algorithm.
Our dual algorithm explicitly considers the dual problem (\ref{feasibility:dual}) and rescales the \emph{row space} of $A$, unlike the previous algorithms.
We already know that the primal algorithm runs faster on the primal infeasible instances.
Thus we expect the dual algorithm runs faster on the dual infeasible instances (i.e., primal feasible instances).
As expected, our dual algorithm does work.
Therefore, in Section \ref{sec:ppdp}, we integrate our primal algorithm and dual algorithm to obtain a quite \emph{balanced} primal-dual algorithm and its performance is also remarkably better than the previous algorithms.

Similar to our primal algorithm, the dual algorithm also consists of the basic procedure BP and the main algorithm MA.
The details of BP and MA are provided in Section \ref{subsec:bpdual} and Section \ref{subsec:madual} respectively.
Similar to \eqref{feasibility:normalized primal}, we consider the normalized version of (\ref{feasibility:dual}) due to the homogeneity:
\begin{equation}
\label{eq:bounded dual}
\mathrm{find }\; u\in \R^m \; \mathrm{subject} \; \mathrm{to} \; x=A^T u> 0, x\in (0,1]^n.
\end{equation}

\subsection{Basic Procedure for the Dual Problem}
\label{subsec:bpdual}

The basic procedure for the dual problem is described in Algorithm \ref{code:dual Dadush}.
\begin{algorithm}[!htb]
	\caption{Basic Procedure for the Dual Problem}
	\label{code:dual Dadush}
	\begin{algorithmic}[1]
        \REQUIRE
          $P_A$
        \ENSURE
            $y, z, J, \case.$
        \STATE $\ths=1/2n^{3/2}, c\in (0,2), \case=0$
        \STATE  $y = \1/n, z=P_Ay, v=y-z=Q_Ay$
		\WHILE {$\case=0$}
		\IF {$v> 0$}
            \STATE $\case=2$ (dual feasible); \rt
		\ELSIF {$z\geq 0$}
		    \STATE $\case=1$ (primal feasible); \rt
		\ELSE
		    \STATE find $K=\{k:\langle z, e_k\rangle\leq 0\}$
		    \STATE $p_K=P_A \sum_{k\in K}e_k$
		    \STATE $\alpha=\langle \frac{p_K}{\|p_K\|_2}, z\rangle $
		    \IF {$\alpha \leq -\ths$}
		        \STATE $y=y-c(\frac{\alpha}{\|p_K\|_2}\sum_{k\in K} e_k)$
                \STATE $z=z-c(\frac{\alpha}{\|p_K\|_2} \sum_{k\in K} p_k)$
                \STATE $v=y-z$
		    \ELSE
		        \STATE find a nonempty set $J$ such that $J\subseteq \{j: \bound'_j(y)\leq \frac{1}{2} \}$ (a cut); \rt
		    \ENDIF
		\ENDIF
		\ENDWHILE
	\end{algorithmic}
\end{algorithm}

In this basic procedure, either a feasible solution for the primal problem is found, or a dual feasible solution is found, or a cut of the bounded row space is found (which is denoted as $bound'_j(y)$ in Line 17).
Now, we need to provide an upper bound in Lemma \ref{lm:dualcut}, which shows that a cut of the bounded row space can be derived from $z$, instead of $v$ in the case of null space (see Lemma \ref{lm:roosbound}).
\begin{lemma}
\label{lm:dualcut}
	Let $x$ be any feasible solution of \eqref{eq:bounded dual} and $z=P_A y$ for some $y$. Then every non-zero coordinate $z_j$ of $z$ gives rise to an upper bound for $x_j$, according to
	\[
	x_j \leq \bound'_j(y) \triangleq \bm{1}^T \Big[\frac{z}{-z_j}\Big]^+.
	\]
\end{lemma}
The proof of this lemma is deferred to Appendix \ref{app:dualcut}.
According to this lemma, we can obtain the following guaranteed cut in Lemma \ref{lm:cutd}, which is similar to Lemma \ref{lm:cut}.
The proof is almost the same as Lemma \ref{lm:cut} just by replacing Lemma \ref{lm:roosbound} with our Lemma \ref{lm:dualcut}.

\begin{lemma}\label{lm:cutd}
Let $P_A$ be the projection matrix at a given iteration of BP (Algorithm \ref{code:dual Dadush}).
Suppose that $\alpha=\inner{\frac{p_K}{\|p_K\|_2}}{z}> -\ths$, then the set $J=\{j: \bound'_j(y)\leq \frac{1}{2} \}$ is nonempty and every solution $x$ of problem (\ref{feasibility:normalized primal}) satisfies $x_j\leq \frac12$ for all $j\in J$.
\end{lemma}

Same to Algorithm \ref{alg:bp}, for the time complexity $\tbp$ of the dual Algorithm \ref{code:dual Dadush}, we have the following lemma.
\begin{lemma}\label{lm:tbpd}
The time complexity of Algorithm \ref{code:dual Dadush} $\tbp=O(n^3m)$. Concretely, it uses at most $O(n^2)$ iterations and each iteration costs at most $O(mn)$ time.
\end{lemma}
Note that Lemma \ref{lm:tbp} also holds regardless \eqref{eq:bounded dual} is feasible or infeasible.
Now, we want to point out that our dual algorithm can perform much better on the dual infeasible instances (primal feasible instances) under well-condition as we expected and discussed before.
Similar to Lemma \ref{lm:infeasible}, we have the following lemma for the dual algorithm.
\begin{lemma}\label{lm:infeasibled}
If \eqref{eq:bounded dual} is infeasible,
the time complexity of Algorithm \ref{alg:bp} $\tbp=O(n^2m/\rho(P_A))$, where $\rho(P_A)$ is the condition number defined in (\ref{eq:rho}). 
In particular, $\rho(P_A)$ equals to $1/\sqrt{n}$ under well-condition (e.g., $A=(I,-I)$, where $I$ is an identity matrix), then $\tbp=O(n^{2.5}m)$ if problem \eqref{eq:bounded dual} is infeasible.
\end{lemma}
Note that it is easy to see that \eqref{feasibility:normalized primal} is feasible (i.e., \eqref{eq:bounded dual} is infeasible) if $A=(I,-I)$.

Besides,
when the dual problem \eqref{eq:bounded dual} is feasible,  we can also utilize the geometry of the problem to bound the iteration complexity instead of Lemma \ref{lm:tbpd}.
Consider the following kind of condition number of the set $Im(A)\bigcap [0,1]^n$:
\[
\delta_{\infty}(Im(A)\bigcap [0,1]^n)\triangleq \max_{x}\{\prod_{i} x_i : x\in Im(A)\bigcap [0,1]^n\}.
\]

As each rescaling in the basic procedure will at least enlarge the value of $\delta_{\infty}(Im(A)\bigcap [0,1]^n)$ by two times, and the largest possible value of $\delta_{\infty}(Im(B)\bigcap [0,1]^n)$ for all matrices $B$ is 1, it takes at most $-\log_2 \delta_{\infty}(Im(A)\bigcap [0,1]^n)$ basic procedures before getting a feasible solution. This means that the iteration complexity of the whole algorithm is $O(n^2 \log \frac{1}{\delta_{\infty}(Im(A)\bigcap [0,1]^n)})$.

\subsection{Main Algorithm for the Dual Problem}\label{subsec:madual}

The main algorithm for the dual problem is described in Algorithm \ref{code:main procedure of dual Dadush}.
Particularly, we rescale the row space of $A$ in Line 6.

\begin{algorithm}[!htb]
	\caption{Main Algorithm for the Dual Problem}
	\label{code:main procedure of dual Dadush}
	\begin{algorithmic}[1]
        \REQUIRE
            $A\in\R^{m\times n}, d=\1, \tau=2^{-L}, \case=0$.
		\WHILE{$\case=0$}
            \STATE $P_A=I-A^T (A A^T)^{\dag}A$
            \STATE $(y,z,J,\case)\leftarrow$ Basic Procedure for Dual Problem$(P_A)$
		\IF {$\case=0$}
		    \STATE $d_J=2d_J$
            \STATE $A_J=2A_J$
		    \STATE $H=\{i:d_i\geq 2^L\}$
            \STATE $d_H=0$
		    \STATE $A_H=0$
		\ENDIF
		\ENDWHILE
        \IF{$\case=1$}
            \STATE $D=\diag(d)$
            \STATE $x=Dz$
        \ENDIF
	\end{algorithmic}
\end{algorithm}

\newpage
Now, we have the following theorem for our dual algorithm.
The proof is provided in Appendix \ref{app:thmdual}.
\begin{theorem}\label{thm:dual}
The time complexity of the dual algorithm is $O(n^4mL)$.
\end{theorem}

\section{Our Primal-Dual $\ppdp$ Algorithm}
\label{sec:ppdp}
In this section, we propose a new polynomial primal-dual projection algorithm (called $\ppdp$) to take advantages of our primal algorithm and dual algorithm.
Similarly, the $\ppdp$ algorithm also consists of two procedures (MA and BP).
Intuitively, the BP solves problems (\ref{feasibility:primal}) and (\ref{feasibility:dual}) simultaneously.
Recall that the primal algorithm runs faster on the infeasible instances and the dual algorithm runs faster on the feasible instances (see Table \ref{table:comparison}).
The MA rescales the matrix (row space or null space) according to the output of BP.
The MA and BP are formally described in Algorithms  \ref{code:main procedure of primal-dual} and \ref{code:primal dual} respectively. The details are deferred to Appendix \ref{app:algo}.
Thus, we have the following theorem.

\begin{theorem}\label{thm:primal_dual}
The time complexity of our primal-dual $\ppdp$ algorithm is $O(n^4mL)$.
\end{theorem}

Note that the final output of our $\ppdp$ algorithm is a feasible solution for either (\ref{feasibility:primal}) or  (\ref{feasibility:dual}).
Obviously, the algorithm will stop whenever it finds a solution of \eqref{feasibility:primal} or \eqref{feasibility:dual}, thus the time complexity of our $\ppdp$ algorithm follows easily from Theorems \ref{thm:primal} and \ref{thm:dual}.

\section{Experiments}
\label{sec:exp}
In this section, we compare the performance of our algorithms with Roos' algorithm \citep{roos2018improved} and Gurobi (one of the fastest solvers nowadays).
We conduct the experiments on the randomly generated matrices.
Concretely, we generate $100$ integer matrices $A$ of size $625\times 1250$, with each entry uniformly randomly generated in the interval $[-100,100]$. The parameter $c\in (0,2)$ is the step-size which is a new practical term introduced in this work.
The average running time of these algorithms are listed in Table \ref{table:comparison}.

\begin{table}[!htb]
  \caption{Running time (sec.) of algorithms wrt. (\ref{feasibility:primal}) is feasible or infeasible}
  \vspace{1mm}
 \label{table:comparison}
 \centering
 \begin{tabular}{ccc}
  \toprule
  Algorithms  & 	feasible instances & 	infeasible instances \\
  \midrule
  Gurobi (a fast optimization solver)  & 3.08  &1.58 \\
  Roos's algorithm \citep{roos2018improved} 
  &  10.75          & 0.83     \\
  Our primal algorithm ($c=1.8$) & 9.93     &  \red{0.48}  \\
  Our dual algorithm ($c=1.8$) & \red{0.35}          &  4.57    \\
  Our $\ppdp$ algorithm ($c=1.8$) & \red{0.60}          & \red{0.58}     \\
  \bottomrule
 \end{tabular}
\end{table}

Table \ref{table:comparison} validates that our new primal-dual $\ppdp$ algorithm is quite balanced on the feasible and infeasible instances due to the integration of our primal and dual algorithm.
Moreover, it shows that our $\ppdp$ algorithm can be a practical option for linear programming since it runs remarkably faster than the fast optimization solver Gurobi.

\section{Conclusion}
\label{sec:con}
In this paper, we try to theoretically explain why the Chubanov-type projection algorithms usually run much faster on the primal infeasible instances.
Furthermore, to address this unbalanced issue, we provide a new fast polynomial primal-dual projection algorithm (called $\ppdp$) by integrating our primal algorithm (which runs faster on the primal infeasible instances) and our dual algorithm (which runs faster on the primal feasible instances).
As a start, we believe more improvements (e.g., the amortized analysis speedup) can be made for the Chubanov-type projection algorithms both theoretically and practically.

\section*{Acknowledgments}
We would like to thank Jian Li (Tsinghua University), Yuanxi Dai (Tsinghua University) and Rong Ge (Duke University) for useful discussions.

\bibliographystyle{plainnat}
\bibliography{bib}

\newpage
\appendix

\section{Details of $\ppdp$ Algorithm}
\label{app:algo}

In this appendix, we describe the details of our $\ppdp$ Algorithm.
Concretely, the main procedure is an integration of Algorithms \ref{alg:ma} and \ref{code:main procedure of dual Dadush}, which is formally described in Algorithm \ref{code:main procedure of primal-dual}.
The basic procedure is an integration of Algorithms \ref{alg:bp} and \ref{code:dual Dadush}, which is formally described in Algorithm \ref{code:primal dual}.

\begin{algorithm}[!h]
	\caption{Main Algorithm for Primal-Dual Problem}
	\label{code:main procedure of primal-dual}
	\begin{algorithmic}[1]
		\REQUIRE $A, d=\1, y=\1/n, \tau=2^{-L}, \case=0$.
        \STATE $A1=A2=A$
        \STATE $d1=d2=d$
        \STATE $y_1=y_2=y$
		\WHILE {case$\leq 0$}
        \STATE $P_{A1}=I-A1^T (A1 A1^T)^{\dag}A1$
		\STATE $Q_{A2}=A2^T (A2 A2^T)^{\dag}A2$
		\STATE $(y_1,y_2,J_1,J_2,\case) \leftarrow$ Basic Procedure for Primal-Dual Problem$(P_{A1},Q_{A2},y_1,y_2)$
		\IF  {$\case=-2$}
        \STATE $d2_{J_2}=d2_{J_2}/2$
        \STATE $H=\{i:d2_i\leq \tau\}$
		\STATE $d2_H=0$
		\STATE $A2_{J_2}=A2_{J_2}/2$
		\STATE $A2=A2_{\overline{H}}$
		\STATE $r=size(Q_{A2})$
		\STATE $y_2=\1/r$

		\ELSIF {$\case=-1$}
		    \STATE $d1_{J_1}=2d1_{J_1}$
            \STATE $A1_{J_1}=2A1_{J_1}$
            \STATE $y_1=y$
        \ENDIF
        \ENDWHILE
		\IF {$\case=1$}
		 \STATE $D=\diag(d1)$
		 \STATE $x=D z_1$
		
        \ELSIF {$\case=2$}
            \STATE $D=\diag(d1)$
            \STATE $x=Dv_1$
        \ELSIF {$\case=3$}
           \STATE $d=d2_{\overline{H}}$
           \STATE $D=\diag(d)$
           \STATE Define $x$ as $x_{\overline{H}}=Dz_2, x_H=0$
        \ELSIF {$\case=4$}
            \STATE $D=\diag(d2)$
            \STATE $x=Dv_2$
		\ENDIF
		
	\end{algorithmic}
\end{algorithm}

\begin{algorithm}[!h]
	\caption{Basic Procedure for Primal-Dual Problem($P_{A1},Q_{A2},y_1,y_2$)}
	\label{code:primal dual}
	\begin{algorithmic}[1]
        \REQUIRE
        $P_{A1},Q_{A2},y_1,y_2$
        \ENSURE
        $y_1, y_2, J_1, J_2, \case.$
        \STATE  $\ths_1=1/2n^{3/2}, r=size(Q_{A2}), \ths_2=1/2r^{3/2}, c\in (0,2), \case=0$
        \STATE $z_1=P_{A1}y_1, v_1=y_1-z_1, v_2=Q_{A2}y_2, z_2=y_2-v_2$
		\WHILE {$\case=0$}
		\IF {$v_1> 0$}
		    \STATE $\case=2$ ($y_1$ is dual feasible); \rt
        \ELSIF {$z_1> 0$}
            \STATE $\case=1$ ($y_1$ is primal feasible); \rt
        \ELSIF {$v_2> 0$ and $r==n$}
            \STATE $\case=4$ ($y_2$ is dual feasible); \rt
		\ELSIF {$z_2> 0$}
		    \STATE $\case=3$ ($y_2$ is primal feasible); \rt
		\ELSE
		    \STATE find $K_1=\{k:\langle z_1, e_k\rangle\leq 0\}$
		    \STATE $p_{K_1}=P_{A1} \sum_{k\in K_1}e_k$
		    \STATE $\alpha_1=\langle \frac{p_{K_1}}{\|p_{K_1}\|_2}, z_1\rangle $
            \STATE find $K_2=\{k:\langle v_2, e_k\rangle\leq 0\}$
		    \STATE $q_{K_2}=Q_{A2} \sum_{k\in K_2}e_k$
		    \STATE $\alpha_2=\langle \frac{q_{K_2}}{\|q_{K_2}\|_2}, v_2\rangle $
		    \IF {$\alpha_1 \leq -\ths_1$}
		        \STATE $y_1=y_1-c(\frac{\alpha_1}{\|p_{K_1}\|_2}\sum_{k\in K_1} e_k)$
		        \STATE $z_1=z_1-c(\frac{\alpha}{\|p_{K_1}\|_2} \sum_{k\in K_1} p_k)$
                \STATE $v_1=y_1-z_1$
		    \ELSE
		        \STATE find a nonempty set $J_1$ such that $J_1\subseteq \{j: \bound'_j(y_1)\leq \frac{1}{2} \}$
                \STATE $\case=-1$; (the rescaling should be in the row space of $A1$); \rt
		    \ENDIF
            \IF {$\alpha_2 \leq -\ths_2$}
		        \STATE $y_2=y_2-c(\frac{\alpha_2}{\|q_{K_2}\|_2}\sum_{k\in K_2} e_k)$
		        \STATE $v_2=v_2-c(\frac{\alpha}{\|q_{K_2}\|_2} \sum_{k\in K_2} q_k)$
                \STATE $z_2=y_2-v_2$
		    \ELSE
		        \STATE find a nonempty set $J_2$ such that $J_2\subseteq \{j: \bound_j(y_2)\leq \frac{1}{2} \}$
                \STATE $\case=-2$; (the rescaling should be in the null space of $A2$); \rt
		    \ENDIF
		\ENDIF
		\ENDWHILE
	\end{algorithmic}
\end{algorithm}

\clearpage
\section{Missing Proofs}
\label{app:pf}

In this appendix, we provide all the proofs for Theorems \ref{thm:primal}--\ref{thm:dual} and Lemmas \ref{lm:cut}--\ref{lm:dualcut}.
\subsection{Proof of Lemma \ref{lm:cut}}
\label{app:lmcut}

As the basic procedure has not terminated, $z=y-v$ cannot be primal feasible. Initially $y=\1/n$.
For each iteration, the operation
\[
y=y-c\big(\frac{\alpha}{\|q_K\|_2}\sum_{k\in K} e_k\big)
\]
will only increase some components of the vector $y$ by $-c\frac{\alpha}{\|q_K\|_2}\geq c\frac{\ths}{\|q_K\|_2}$. Thus each component of $y$ is at least $1/n$ during the whole procedure， implying that there exists an index $j$ with $v_j\geq 1/n$. Otherwise $z=y-v$ will be primal feasible. For this specific $j$, we have
\begin{equation*}
\begin{aligned}
\bound_j(y)&=\frac{-\sum_{i=1}^{n} v_i^{-}}{v_j}=\frac{-\sum_{k\in K} \langle e_k, v  \rangle}{v_j}\\
&=\frac{-\sum_{k\in K} \langle Q_A e_k, v  \rangle}{v_j}=\frac{-\langle q_K, v  \rangle}{v_j}\\
&\leq n\|q_K\|_2\ths\leq\frac{1}{2},
\end{aligned}
\end{equation*}
where the last inequality follows from
\[
\|q_K\|_2\leq \|\sum_{k\in K} e_k\|\leq \sqrt{n}.  \]
$\QED$

\subsection{Proof of Lemma \ref{lm:tbp}}\label{app:tbp}
When $\alpha=\langle \frac{q_K}{\|q_K\|_2}, v\rangle \leq -\ths$, the decrease of $\|v\|_2$ in each iteration has a lower bound:
\[
\|v-c\frac{\alpha}{\|q_K\|_2} q_K\|_2^2=\|v\|_2^2-(2c-c^2)\alpha^2\leq \|v\|_2^2-(2c-c^2)\frac{1}{4n^3}.
\]
Initially, $\|v_0\|_2^2\leq \|y_0\|_2^2=\frac{1}{n}$. After $t$ iterations, $\|v_t\|_2^2\leq \frac{1}{n}-\frac{(2c-c^2)t}{4n^3}$. So it takes at most $O(n^2)$ iterations to obtain a vector $v_t$ with $\|v_t\|_2\leq \frac{1}{n}$, in which case a primal feasible solution $z=y-v$ can be obtained. This means that each basic procedure takes at most $O(n^2)$ iterations before stopping.

Now it remains to bound the time complexity in each iteration. In each basic procedure iteration, we find all indices $k$ such that $v_k\leq 0$ and do the calculation $\sum_{k\in K} q_k $.
In the worst case, $|K|$ can be $O(n)$, thus $O(n^2)$ arithmetic operations are needed to calculate this summation. However, there is another way to do this. Recall that $\rank(Q_A)=\rank(A)=m$, thus the number of basis vectors of the $n$ rows are exactly $m$, while the other rows can be represented by a weighted summation of these $m$ basis vectors. This means that we only have to do the naive summation for these $m$ rows, while the other $n-m$ rows can be obtained by doing the weighted summation of these $m$ elements， which cost $O(mn)+O((n-m)m)=O(mn)$ arithmetic operations. Note that the $m$ basis vectors can be computed by the Singular Value Decomposition (SVD). Then the weights of the other $n-m$ rows can be obtained by computing the inverse of an $m\times m$ matrix and multiplying this inverse matrix to the row vectors. These steps cost $O(mn^2+m^3+m^2(n-m))=O(mn^2)$ operations and only needs to be done once at the beginning of the basic procedure.
$\QED$

\subsection{Proof of Lemma \ref{lm:infeasible}}\label{app:infeasible}

When (\ref{feasibility:normalized primal}) is infeasible, i.e., $\rho(Q_A)>0$. Denoting the vector $w$ as the center which achieves the value $\rho(Q_A)$, we can check the closeness between $v_t$ and $w$ in each iteration:
\[
\langle v_{t+1},w  \rangle =\langle v_{t},w  \rangle -c\langle v_{t},\hat{q_K}  \rangle \langle w,\hat{q_K} \rangle \geq \langle v_{t},w  \rangle+ \frac{c\rho(Q_A)}{2n\sqrt{n}}.
\]
On the other hand, as the norm $\|v_t\|\leq \sqrt{\frac{1}{n}-\frac{(2c-c^2)t}{n^3}}$ before the basic procedure stops, we should have
\[
\rho(Q_A)+\frac{ct\rho(Q_A)}{2n\sqrt{n}}\leq \langle v_{t},w  \rangle \leq \|v_{t}\|\|w\|\leq \sqrt{\frac{1}{n}-\frac{(2c-c^2)t}{n^3}}.
\]
This implies that the number of the iterations in the basic procedure is $t=O(\min\{\frac{n}{\rho(Q_A)},n^2\})$ when the primal problem (\ref{feasibility:normalized primal}) is infeasible.
This proof is finished by combining the result of Lemma \ref{lm:tbp}, i.e., each iteration costs $O(mn)$ time.
$\QED$

\subsection{Proof of Theorem \ref{thm:primal}}\label{app:thmprimal}
According to a classic result of \citep{khachian1979polynomial}, there exists a positive number $\tau$ (satisfying $1/\tau=O(2^L)$) such that the positive coordinates of the basic feasible solutions of problem \eqref{feasibility:normalized primal} are bounded below by $\tau$.  When the coordinate has been rescaled for more than $\tau$ times, the value of this coordinate in all the solutions must be $0$. As a result, the corresponding columns of $A$ can be omitted.

According to Lemma \ref{lm:cut}, the cut $J$ is nonempty. It means each iteration of MA can rescale at least one coordinate of the feasible solutions by $1/2$. Thus the number of rescaling operations can be bounded by $nL$, i.e., the number of iterations in MA $\tma=O(nL)$. The proof is finished by combining this with Lemma \ref{lm:tbp}.
$\QED$

\subsection{Proof of Lemma \ref{lm:dualcut}}\label{app:dualcut}

Since $x\in \mathcal{R}_A$, we have $\langle x,z  \rangle=0$.
Thus we consider the following two cases.
For $z_j< 0$, we have
\[
-z_j x_j=\sum_{i\neq j} z_i x_i \leq \sum_{i, z_i>0}z_i x_i\leq \sum_{i, z_i>0}z_i=\bm{1}^T z^+.
\]
On the other hand, for $z_j> 0$, we have
\[
z_j x_j=-\sum_{i\neq j} z_i x_i \leq \sum_{i, z_i<0}-z_i x_i\leq \sum_{i, z_i<0}-z_i=-\bm{1}^T z^-.
\]
$\QED$

\subsection{Proof of Theorem \ref{thm:dual}}\label{app:thmdual}
First, we note that (\ref{feasibility:dual}) is feasible if and only if the problem
\begin{equation}\label{eq:dddual}
\begin{aligned}
&P_A x= 0,\quad
x> 0\quad
\end{aligned}
\end{equation}
is feasible. The reason is that (\ref{feasibility:dual}) is feasible if and only if there is a positive vector $x=A^Tu>0$ in the row space of $A$.
Further, this system \eqref{eq:dddual} is feasible if and only if the following normalized system
\begin{equation}
\label{eq:normalized dual}
\begin{aligned}
P_A x= 0,\quad
x\in (0,1]^n
\end{aligned}
\end{equation}
is feasible.
In problem \eqref{eq:normalized dual}, the bit length of $P_A$ can be the same as $A$, i.e.,  $O(L)$.
Now, the remaining proof is the same as that for Theorem \ref{thm:primal}.
According to the classic result of \citep{khachian1979polynomial},  there exists a positive number $\tau$ (satisfying $1/\tau=O(2^L)$) such that the positive coordinates of the basic feasible solutions of problem \eqref{eq:normalized dual} are bounded below by $\tau$.
Also according to Lemma \ref{lm:cutd}, the cut $J$ is nonempty. It means each iteration of MA can rescale at least one coordinate of the feasible solutions by $1/2$. Thus the number of rescaling operations can be bounded by $nL$, i.e., the number of iterations in MA $\tma=O(nL)$. The proof is finished by combining this with Lemma \ref{lm:tbpd}.
$\QED$

\end{document}